\documentclass{amsart}

\usepackage{amssymb,amscd,amsthm,latexsym}
\usepackage[all]{xy}
\usepackage{color}

\def\cleft{\hbox{[\kern-.16em\hbox{[}}}
\def\cright{\hbox{]\kern-.16em\hbox{]}}}

\def\CC{ \mathbb{C} }

\def\EE{ \mathbb{E} }

\newdir{ >}{{}*!/-8pt/@{>}}

\def\into{ \rightarrowtail }
\def\onto{ \twoheadrightarrow }
\def\splito{\rightleftarrows }

\def\pullback{
 \ar@{-}[]+R+<6pt,-1pt>;[]+RD+<6pt,-6pt>
 \ar@{-}[]+D+<1pt,-6pt>;[]+RD+<6pt,-6pt>}

\newcommand{\lhdi}{\lhd^{-1}}
\def\cleft{\hbox{[\kern-.16em\hbox{[}}}
\def\cright{\hbox{]\kern-.16em\hbox{]}}}

\theoremstyle{plain}
\newtheorem{defi}[subsection]{Definition}
\newtheorem{theo}[subsection]{Theorem}
\newtheorem{lma}[subsection]{Lemma}
\newtheorem{prop}[subsection]{Proposition}
\newtheorem{coro}[subsection]{Corollary}

\theoremstyle{remark}

\newtheorem{exms}[subsection]{Examples}

\begin{document}

\title[The $3 \times 3$ lemma in the $\Sigma$-Mal'tsev and $\Sigma$-protomodular settings]{The $3 \times 3$ lemma in the $\Sigma$-Mal'tsev
and $\Sigma$-protomodular settings. Applications to monoids and
quandles}

\author{Dominique Bourn}
\address[Dominique Bourn]{LMPA, Universit\'{e} du Littoral, 50 Rue Ferdinand Buisson, BP 699, 62228 Calais, France}
\thanks{}
\email{bourn@lmpa.univ-littoral.fr}

\author{Andrea Montoli}
\address[Andrea Montoli]{Dipartimento di Matematica ``Federigo Enriques'', Universit\`{a} degli Studi di Milano, Via Saldini 50, 20133 Milano, Italy}
\thanks{} \email{andrea.montoli@unimi.it}



\begin{abstract}
We investigate what is remaining of the $3 \times 3$ lemma and of
the denormalized $3 \times 3$ lemma, respectively valid in a
pointed protomodular and in a Mal'tsev category, in the context of
partial pointed protomodular and partial Mal'tsev categories,
relatively to a class $\Sigma$ of points. The results apply, among
other structures, to monoids, semirings and quandles.
\end{abstract}

\maketitle

\section{Introduction}
The $3\times 3$ lemma is a classical tool in homological algebra,
with several applications. It holds for many algebraic structures,
including groups. As shown in \cite{Bourn 3x3}, the lemma is valid
in any pointed regular protomodular \cite{Bourn protomod}
category: given any commutative diagram as the one on the
left-hand side here below, where the three columns and the middle
row are exact sequences, the upper row is exact if and only if the
lower one is.
\begin{equation} \label{diagrams norm and denorm 3x3}
\xymatrix@=25pt{ {K[\phi]\;} \ar@{>->}[r]^{K(k_x)}\ar@{
>->}[d]_{k_{\phi}}  & K[f] \ar[r]^{K(x)} \ar@{ >->}[d]^{k_{f}} &
K[f'] \ar@{ >->}[d]^{k_{f'}} && R[\phi]
\ar@<-1,ex>[d]_{d_0^{\phi}} \ar@<+1,ex>[d]^{d_1^{\phi}}
\ar@<-1,ex>[r]_{R(d_1^x)}\ar@<+1,ex>[r]^{R(d_0^x)}
    & R[f] \ar@<-1,ex>[d]_{d_0^f}\ar@<+1,ex>[d]^{d_1^f} \ar[l] \ar[r]^{R(x)} & R[f'] \ar@<-1,ex>[d]_{d_0^{f'}}\ar@<+1,ex>[d]^{d_1^{f'}}\\
    {K[x]\;} \ar@{->>}[d]_{\phi} \ar@{>->}[r]^{k_x}  &  X  \ar@{->>}[r]^{x} \ar@{->>}[d]^{f} & X'  \ar@{->>}[d]^{f'} && R[x] \ar@{->>}[d]_{\phi}\ar[u]
    \ar@<-1,ex>[r]_{d_0^x}\ar@<+1,ex>[r]^{d_1^x} & X \ar@{->>}[d]^{f}\ar[u] \ar[l] \ar@{->>}[r]^{x} & X' \ar@{->>}[d]^{f'}\ar[u]\\
    U \ar[r]_u &  Y   \ar@{->>}[r]_y  & Y'    &&    W \ar@<-1,ex>[r]_{y_0} \ar@<+1,ex>[r]^{y_1}  & Y  \ar[l] \ar@{->>}[r]_{y} & Y'
}
\end{equation}
A \emph{denormalized} version of the $3\times 3$ lemma was proved
in \cite{Bourn denormalized 3x3} to be valid in any (not
necessarily pointed) regular Mal'tsev category: given any
commutative diagram as the one on the right-hand side, where the
three columns and the middle row are exact forks, the upper row is
exact if and only if the lower one is. Later, the denormalized $3
\times 3$ lemma was extended to the context of regular Goursat
categories in \cite{Lack}. A comparison between the normalized and
the denormalized versions of the $3 \times 3$ lemma was described
in \cite{GJR 3x3 lemma}.

In the recent paper \cite{MMS nine lemma} it was shown that a
particular version of the $3 \times 3$ lemma holds also in the
category $Mon$ of monoids, which is not protomodular nor Mal'tsev:
the $3 \times3$ lemma holds for monoids when we replace exact
sequences by \emph{special Schreier} exact sequences. This notion,
introduced in \cite{Schreier book}, originated from the notion of
\emph{Schreier split epimorphism} \cite{MartinsMontoliSobral mon w
op}: these are the split epimorphisms that correspond to classical
monoid actions. An action of a monoid $B$ on a monoid $X$ is a
monoid homomorphism $B \to End(X)$. The $3 \times 3$ lemma for
special Schreier exact sequences allowed to give in \cite{MMS nine
lemma} a description of a Baer sum construction of special
Schreier exact extensions with abelian kernel, obtained thanks to
a \emph{push forward} construction, analogous to the classical one
for group extensions.

In order to understand and to describe categorically the
(co)homological features of Schreier monoid extensions, the
notions of pointed $\Sigma$-protomodular \cite{BMMS
S-protomodular} and $\Sigma$-Mal'tsev category \cite{Bourn sigma
Maltsev} have been introduced, with respect to a class $\Sigma$ of
\emph{points}, i.e. split epimorphisms with a fixed section. The
main examples of the first notion are the categories of monoids,
of semirings (see \cite{BMMS S-protomodular}) and, more generally,
any J\'{o}nsson-Tarski variety \cite{MartinsMontoli}. An
interesting example of a $\Sigma$-Mal'tsev category (which is not
$\Sigma$-protomodular) is the category of quandles \cite{Bourn
quandles}. As for the case of monoids, in all these contexts the
class $\Sigma l$ of $\Sigma$-special maps (see Definition
\ref{special}) appeared to be very discriminating.

The aim of the present paper is to investigate what is remaining
of the normalized and of the denormalized $3 \times 3$ lemma,
respectively in the abstract context of pointed
$\Sigma$-protomodular and of $\Sigma$-Mal'tsev categories. Namely,
we are interested in the description of the conditions under
which, given any (normalized or denormalized) diagram as
\eqref{diagrams norm and denorm 3x3} above, where the three
columns and the middle row are exact (or $\Sigma$-special), the
upper row is exact (or $\Sigma$-special) when the lower one is
(the so-called \emph{upper $3 \times 3$ lemma}), and conversely
(the so-called \emph{lower $3 \times 3$ lemma}).

In these relative contexts, a curious phenomenon appears, in
contrast to the ``absolute'' case of protomodular and Mal'tsev
categories (that are $\Sigma$-protomodular and $\Sigma$-Mal'tsev,
respectively, for the class $\Sigma$ of all points). In fact, in
the absolute contexts, the upper and the lower $3 \times 3$ lemmas
are equivalent, both in the normalized and in the denormalized
case (see \cite{ZJanelidze} and \cite{GranRodelo}). In the
relative contexts, this equivalence is no more true (see Theorem
\ref{th2}, Proposition \ref{p22} and Theorem \ref{cuicui} below).
This shows an unexpected asymmetry between the two parts of the $3
\times 3$ lemma.

The paper is organized as follows. In Section \ref{section
Sigma-Mal'tsev categories} we recall the notion of
$\Sigma$-Mal'tsev category and we give several examples. In
Section \ref{section regular context} we obtain some properties of
regular $\Sigma$-Mal'tsev categories that are used in Section
\ref{section denormalized 3x3 lemma}, where we describe the
versions of the denormalized $3 \times 3$ lemma that are valid in
regular $\Sigma$-Mal'tsev categories. In Section \ref{section
Sigma-protomodular categories} we recall the definition and the
main properties of $\Sigma$-protomodular categories. In Section
\ref{section aspects 3x3 lemma} we describe the versions of the
(normalized) $3 \times 3$ lemma that hold in a regular
$\Sigma$-protomodular category. In Section \ref{section Baer sums}
we give an interpretation of the Baer sum construction in
Barr-exact $\Sigma$-Mal'tsev categories.

\section{$\Sigma$-Mal'tsev categories} \label{section Sigma-Mal'tsev categories}

\subsection{The fibration of points}

Throughout the paper, all the categories we consider will be
finitely complete. A (generalized) \emph{point} in a category
$\EE$ is a pair $(f, s)$ of morphisms such that $fs = 1$; in other
terms, $f$ is a split epimorphism with fixed section $s$. The
category $Pt\EE$ is the category whose objects are the points in
$\EE$ and whose morphisms are the pairs $(y,x)$ of morphisms which
make a square as below commutative, both downward and upward:
\[ \xymatrix{ X' \ar@<-2pt>[d]_{f'} \ar[r]^{x} & X
\ar@<-2pt>[d]_{f} \\
Y' \ar@<-2pt>[u]_{s'} \ar[r]_y & Y. \ar@<-2pt>[u]_s } \]

The codomain functor $\P_{\EE} \colon Pt\EE \to \EE$ is a
fibration whose cartesian maps are the pullbacks of split
epimorphims; it is called \emph{the fibration of points}
\cite{Bourn protomod}. As usual, we denote by $\EE^2$ the category
whose objects are the arrows in $\EE$ and whose morphisms are the
commutative squares.

Let $\Sigma$ be a pullback stable class of points in a category
$\EE$. We denote by $\Sigma Pt\EE$ the full subcategory of $Pt\EE$
whose objects are the points in $\Sigma$; the restriction of
$\P_{\EE}$ to the class $\Sigma$ determines a subfibration of the
fibration of points:
$$\xymatrix{
    \Sigma Pt\EE\,\, \ar@{>->}[r]^j \ar[d]_{\P^{\Sigma}_{\EE}} & Pt\EE  \ar[d]^{\P_{\EE}}  \\
    \EE \ar@{=}[r]  & \EE.}
$$

Recall from \cite{BMMS S-protomodular} the following:
\begin{defi}\label{special}
A reflexive relation on an object $X$:
$$\xymatrix{ R \ar@<-6pt>[r]_{d_1^R} \ar@<6pt>[r]^{d_0^R} & X \ar[l]|{s_0^R} }$$
such that the pair $(d_0^R,s_0^R)$ is a point in $\Sigma$ is
called a \emph{$\Sigma$-relation}. A morphism $f \colon X \to Y$
is said to be \emph{$\Sigma$-special} when its kernel equivalence
relation $R[f]$ is a $\Sigma$-equivalence relation.
\end{defi}
We denote by $\Sigma l$ the class of $\Sigma$-special morphisms.
An object $X$ is said to be \emph{$\Sigma$-special} when the
terminal map $\tau_X \colon X \to 1$ is $\Sigma$-special. {Observe
that, if} a point $(f,s)$ is in $\Sigma$, the morphism $f$ is not
necessarily $\Sigma$-special. However, when a $\Sigma$-special
morphism $f$ is split by $s$, the pair $(f,s)$ is in $\Sigma$ (see
\cite{BMMS S-protomodular} for more details). Finally, we denote
by $\Sigma l\EE^2$ the full subcategory of $\EE^2$ whose objects
are the $\Sigma$-special morphisms, and by $\Sigma l_Y\EE$ its
subcategory whose morphisms have $1_Y$ as lower horizontal map.

\subsection{$\Sigma$-Mal'tsev categories}

Recall from \cite{Bourn sigma Maltsev} the following:
\begin{defi}
A finitely complete category $\EE$ is said to be a
\emph{$\Sigma$-Mal'tsev} category when, given any pullback of
points with $(f,s)\in \Sigma$:
$$\xymatrix{
    X' \pullback \ar@<2pt>[r]^{\bar g} \ar@<-2pt>[d]_{f'} & X \ar@<2pt>[l]^{\bar t} \ar@<-2pt>[d]_f\\
    Y'   \ar@<2pt>[r]^g \ar@<-2pt>[u]_{s'} & Y, \ar@<2pt>[l]^t \ar@<-2pt>[u]_s }
  $$
  the pair $(s',\bar t)$ is jointly extremally epimorphic.
\end{defi}
The previous definition is equivalent to the following one: given
any commutative square of points, with $(f,s)\in \Sigma$:
$$\xymatrix{
    X''  \ar@<2pt>[r]^{\check g} \ar@<-2pt>[d]_{f''} & X \ar@<2pt>[l]^{\check t} \ar@<-2pt>[d]_f\\
    Y'   \ar@<2pt>[r]^g \ar@<-2pt>[u]_{s''} & Y, \ar@<2pt>[l]^t \ar@<-2pt>[u]_s }
  $$
  the unique induced factorization $\phi \colon X''\to X'$ toward the pullback of $(f,s)$ along $g$ is an extremal epimorphism.

\begin{exms}
\begin{enumerate}
  \item In the category $Mon$ of monoids, a point $(f,s) \colon A\splito B$ is \emph{weakly Schreier}
  \cite{Bourn sigma Maltsev} if, for any $b \in B$, the map $\mu_b \colon \text{Ker}(f) \to f^{-1}(b)$ defined by $\mu_b(x) = x\cdot s(b)$ is surjective.
  The point $(f, s)$ is a \emph{Schreier point} \cite{MartinsMontoliSobral mon w op, BMMS SemForum} if, for any $b \in B$, the map $\mu_b$ is
  bijective. It was observed in \cite{BMMS SemForum} that a point $(f,s)$ is a Schreier point if and only if there exists a unique set-theoretical
  map $q_f \colon A \to \text{Ker}(f)$ such that $a = q_f(a) \cdot sf(a)$ for all $a \in A$. The map $q_f$ is called the \emph{Schreier retraction} of
  $(f,s)$.
  As shown in \cite{Schreier book, Bourn sigma Maltsev}, $Mon$ is a
  $\Sigma$-protomodular category, and consequently a $\Sigma$-Mal'tsev one (see Definition \ref{def1} and Theorem \ref{main} below) for $\Sigma$ either
  the class of Schreier or weakly Schreier points. \\

  \item More generally, let $\CC$ be a J\'{o}nsson-Tarski variety, i.e. a variety (in the sense of universal algebra) whose corresponding theory has
  a unique constant $0$ and a binary term $+$ which satisfy the following
  axiom:
  \[ x + 0 = 0 + x = x. \]
  The notion of a (weakly) Schreier point, given as above for monoids, actually makes sense in every J\'{o}nsson-Tarski variety.
  As shown in \cite{MartinsMontoli}, every J\'{o}nsson-Tarski variety $\CC$ is a $\Sigma$-protomodular (and hence $\Sigma$-Mal'tsev) category
  for the class $\Sigma$ of Schreier points.
  Analogously, it is easy to see that $\CC$ is $\Sigma$-protomodular for the class $\Sigma$ of weakly
  Schreier points, too. In particular, the category $SRng$ of semirings is $\Sigma$-protomodular for both classes. \\

  \item The definition of a Schreier point of monoids, and the
  proof that $Mon$ is a $\Sigma$-protomodular category when
  $\Sigma$ is the class of Schreier points, only make use of
  finite limits. Hence they are invariant under the Yoneda
  embedding. This means that it makes sense to consider
  \emph{internal Schreier points} in every category $Mon \EE$ of
  internal monoids in any category $\EE$, and that $Mon \EE$ is
  $\Sigma$-protomodular w.r.t. the class of such points. The same
  is true if we replace $Mon$ with every J\'{o}nsson-Tarski
  variety. \\

  \item A quandle \cite{Joyce, Matveev} is a set $A$ equipped with two binary operations $\lhd$ and $\lhdi$ such that the following
  identities hold (for all $a,\ b,\ c \in A$):
  \begin{itemize}
  \item[(A1)] $a \lhd a = a = a \lhdi a$ (idempotency); \item[(A2)]
  $(a \lhd b) \lhdi b = a = (a \lhdi b) \lhd b$ (right invertibility); \item[(A3)] $(a \lhd b) \lhd c = (a \lhd c) \lhd
  (b \lhd c)$ and $(a \lhdi b) \lhdi c = (a \lhdi c) \lhdi (b \lhdi
  c)$ (self-distributivity).
  \end{itemize}
  The structure of quandle is of interest in knot theory, since the three axioms above correspond to the Reidemeister moves on
  oriented link diagrams. A homomorphisms of quandles is a map which preserves both $\lhd$ and $\lhdi$. In the category $Qnd$ of
  quandles and quandle homomorphisms, a point $(f,s) \colon A\splito B$ is called \emph{puncturing}
  (resp. \emph{acupuncturing}) if, for every $b \in B$, the map $\mu_b \colon f^{-1}(b) \to f^{-1}(b)$, defined by $\mu_b(a) =
  s(b) \lhd a$, is surjective (resp. bijective). In \cite{Bourn quandles} it was shown that $Qnd$ is a $\Sigma$-Mal'tsev category (which is not
  a $\Sigma$-protomodular one)
  w.r.t. both classes of puncturing and acupuncturing points.\\

  \item Let $Cat_B$ be the fiber above a set $B$ of the fibration $(\;)_0 \colon Cat \to Set$ which associates with every category its
  set of objects. A point $(F, S)$ in $Cat_B$ has \emph{fibrant splittings} \cite{Bourn Mal'tsev reflection} if, for every arrow
  $\varphi$ in the codomain of $F$, the arrow $S(\varphi)$ is cartesian. $Cat_B$ is $\Sigma$-protomodular (and hence $\Sigma$-Mal'tsev) for the class
  $\Sigma$ of points with fibrant splittings. We observe that $Mon$ is $Cat_1$, where $1$ is the one-element set, and that in this
  case the notion of point with fibrant splittings reduces to the notion of Schreier point.
  \end{enumerate}
  \end{exms}

We say that the class $\Sigma$ is \emph{point-congruous} when
$\Sigma Pt\EE$ is closed under finite limits in $Pt\EE$ (which
implies that it contains all the isomorphisms); in this case
$\Sigma l\EE^2$ is closed under finite limits in $\EE^2$. If $\CC$
is a J\'{o}nsson-Tarski variety, and $\Sigma$ is the class of
Schreier points, then $\Sigma$ is a point-congruous class. This is
the case, in particular, for the categories $Mon$ of monoids and
$SRng$ of semirings. Similarly, if $\CC = Cat_B$, the class of
points with fibrant splitting is point-congruous. If $\CC =Qnd$,
the class of acupuncturing points is point-congruous. Recall from
\cite{Bourn sigma Maltsev} the following remarkable:

\begin{theo}
Suppose $\Sigma$ is point-congruous and $\EE$ is a
$\Sigma$-Mal'tsev category, then:
\begin{enumerate}
\item when $g f$ and $g$ are $\Sigma$-special, so is $f$; \item
any subcategory $\Sigma l/Y$ of the slice category $\EE /Y$ is a
Mal'tsev category.
\end{enumerate}
\end{theo}

In particular, if we denote by $\Sigma\EE_{\sharp}=\Sigma l_1\EE$
the full subcategory of $\EE$ whose objects are the
$\Sigma$-special objects, it is a Mal'tsev category \cite{BMMS
S-protomodular, Bourn sigma Maltsev}, called the \emph{Mal'tsev
core} of the point-congruous $\Sigma$-Mal'tsev category $\EE$; any
of its morphisms is $\Sigma$-special. If $\CC = Mon$, $\Sigma
\CC_{\sharp}$ is the category $Gp$ of groups; if $\CC = SRng$,
$\Sigma \CC_{\sharp}$ is the category $Rng$ of (not necessarily
unitary) rings (in both cases, we are considering the class
$\Sigma$ of Schreier points). If $\CC = Qnd$, and $\Sigma$ is the
class of acupuncturing points, $\Sigma \CC_{\sharp}$ is the
category of \emph{latin quandles}, i.e. those quandles whose table
of composition for the operation $\lhd$ is a latin square.

\section{The regular context} \label{section regular context}

In this section we suppose that the ground category $\EE$ is
regular. Then the categories $\EE^2$ and $Pt(\EE)$ are regular as
well, and their regular epimorphisms are the levelwise ones.

\begin{defi}[\cite{regpush}] Let $\EE$ be a regular category.
Consider the following commutative diagram:
\begin{equation} \label{diagram regular pushout}
\xymatrix{ & R[f] \ar@<-4pt>[d]_{d_1^f} \ar@<4pt>[d]^{d_0^f}
\ar@{->>}[r]^{R(x)} & R[f'] \ar@<-4pt>[d]_{d_1^{f'}}
\ar@<4pt>[d]^{d_0^{f'}} \\
R[x] \ar@{->>}[d]_{R(f)} \ar@<4pt>[r]^{d_0^x}
\ar@<-4pt>[r]_{d_1^x} & X \ar@{}[dr]|{(*)} \ar[u] \ar[l]
\ar@{->>}[d]_{f}
 \ar@{->>}[r]^{x} & X' \ar[u] \ar@{->>}[d]^{f'}  \\
R[y] \ar@<4pt>[r]^{d_0^y} \ar@<-4pt>[r]_{d_1^y} & Y
\ar@{->>}[r]_{y} \ar[l] & Y',}
\end{equation}
where we denote by $R[x]$ the kernel pair of $x$, and by $R(f)$
the morphism between the kernel  pairs induced by $f$. The
commutative square $(*)$ of regular epimorphisms is said to be a
\emph{regular pushout} whenever the induced factorization $X \to Y
\times_{Y'} X'$ into the pullback of $f'$ along $y$ is a regular
epimorphism.
\end{defi}
\noindent Any regular pushout is in particular a pushout.
Moreover, the factorizations $R(f)$ and $R(x)$ as in
\eqref{diagram regular pushout} are regular epimorphisms.

\begin{prop}\label{p10}
Let $\EE$ be a regular $\Sigma$-Mal'tsev category. Consider any
regular epimorphism of points as on the right-hand side of the
following diagram:
\[ \xymatrix{ R[x] \ar@<-2pt>[d]_{R(f)}
\ar@<4pt>[r]^{d_0^x}
\ar@<-4pt>[r]_{d_1^x} & X  \ar[l] \ar@<-2pt>[d]_{f} \ar@{->>}[r]^{x} &  X' \ar@<-2pt>[d]_{f'}  \\
R[y] \ar@<-2pt>[u]_{R(s)}  \ar@<4pt>[r]^{d_0^y}
\ar@<-4pt>[r]_{d_1^y} & Y \ar@{->>}[r]_{y} \ar[l]
\ar@<-2pt>[u]_{s} & Y'. \ar@<-2pt>[u]_{s'} } \] When its domain
$(f,s)$ is in $\Sigma$, or when the regular epimorphism $y$ is
$\Sigma$-special, the downward square on the right is a regular
pushout, and the factorization \linebreak $R(x) \colon R[f] \to
R[f']$ is a regular epimorphism.
\end{prop}

\proof Consider the following diagram:
\[ \xymatrix{ R[x] \ar@<-2pt>[d]_{R(f)} \ar[r]^{\psi} & \bar{X}
\ar@<-2pt>[d]_{\bar{f}} \pullback \ar[r]^{\bar{d}_0^y} & X
\ar@<-2pt>[d]_{f} \ar[r]^{\phi} & \bar{X}'
\ar@<-2pt>[d]_{\bar{f}'} \pullback \ar[r]^{\bar{y}} & X' \ar@<-2pt>[d]_{f'}  \\
R[y] \ar@<-2pt>[u]_{R(s)} \ar@{=}[r] & R[y]
\ar@<-2pt>[u]_{\bar{s}} \ar[r]_{d_0^y} & Y \ar@{=}[r]
\ar@<-2pt>[u]_{s} & Y \ar@<-2pt>[u]_{\bar{s}'} \ar[r]_y & Y',
\ar@<-2pt>[u]_{s'} }
\] where $(\bar{f}', \bar{s}')$ is the pullback of $(f', s')$ along $y$,
$(\bar{f}, \bar{s})$ is the pullback of $(f,s)$ along $d_0^y$, and
$\psi$ and $\phi$ are induced by the universal property of the
pullback. Since $\EE$ is a $\Sigma$-Mal'tsev category, $\psi$ is a
regular epimorphism whenever $(f,s)$ is in $\Sigma$, or whenever
$(d_0^y,s_0^y)$ is in $\Sigma$, namely when $y$ is
$\Sigma$-special. Denote by $\delta_1 = (x \bar{d}_0^y, d_1^y
\bar{f}) \colon \bar X \to \bar X'$ the factorization into the
pullback. The morphism $\delta_1$ is the pullback of $x$ along
$\bar y$, and so it is a regular epimorphism. Moreover, we get
$\phi \, d_1^x = \delta_1 \, \psi$, which is a regular
epimorphism, being the composite of two regular epimorphisms.
Accordingly, $\phi$ itself is a regular epimorphism, and the
square in question is a regular pushout.
\endproof

Whence an important consequence about the direct images of
equivalence relations along regular epimorphisms; in general they
are no longer equivalence relations, but only reflexive and
symmetric ones. However we get:

\begin{coro}
Let $\EE$ be a regular $\Sigma$-Mal'tsev category, $f \colon X
\onto Y$ a regular epimorphism, and $R$ an equivalence relation on
X. When $R$ is a $\Sigma$-equivalence relation, or when $f$ is
$\Sigma$-special, the direct image $f(R)$ is itself an equivalence
relation on $Y$.
\end{coro}

\proof Consider the following diagram, where $f(R)$ is given by
the canonical decomposition $R \onto f(R)\into Y\times Y$ of the
map $(f \, d_0^R, f \, d_1^R) \colon R \to Y \times Y$:
$$\xymatrix{
R   \ar@<-6pt>[d]_{d_0^R}  \ar@{->>}[r]^{\check f} \ar@<6pt>[d]^{d_1^R} &  f(R) \ar@<-6pt>[d]_{d_0} \ar@<6pt>[d]^{d_1} \\
X \ar@{->>}[r]_{f} \ar[u]|{s_0^R} & Y. \ar[u]|{s_0} }
$$
According to the previous proposition, under the assumption that
$R$ is a $\Sigma$-equivalence relation or that $f$ is
$\Sigma$-special, the split vertical square indexed by $0$ is a
regular pushout. Accordingly, the factorization $R(\check f)
\colon R[d_0^R]\to R[d_0]$ is a regular epimorphism. Thanks to
Lemma $1.7$ in \cite{Bourn Goursat}, we obtain that $f(R)$ is an
equivalence relation.
\endproof

\begin{prop}\label{p1}
Let $\EE$ be a regular $\Sigma$-Mal'tsev category. Consider a
commutative square of regular epimorphisms as  below, where $f$ is
$\Sigma$-special:
\begin{equation} \label{diagram comm square of reg epis}
\xymatrix{
X  \ar@{->>}[d]_{f}  \ar@{->>}[r]^{x} & X'  \ar@{->>}[d]^{f'}  \\
Y \ar@{->>}[r]_{y} & Y'. }
\end{equation}
Suppose moreover that the factorization $R(x) \colon R[f] \to
R[f']$ is a regular epimorphism. Then the square \eqref{diagram
comm square of reg epis} is a regular pushout and the
factorization $R(f) \colon R[x] \to R[y]$ is a regular
epimorphism. In other words, \eqref{diagram comm square of reg
epis} is a regular pushout if and only if $R(x)$ is a regular
epimorphism.
\end{prop}

\proof Complete Diagram \eqref{diagram comm square of reg epis} by
the kernel equivalence relations of the vertical maps:
$$\xymatrix{
R[f] \ar@<-6pt>[d]_{d_0^f}  \ar@{->>}[r]^{R(x)}
\ar@<6pt>@{.>}[d]^{d_1^f} & R[f'] \ar@<-6pt>[d]_{d_0^{f'}}
\ar@<6pt>@{.>}[d]^{d_1^{f'}} \\
X \ar@{->>}[r]_{x} \ar@{->>}[d]_{f} \ar[u]|{s_0^f} & X'
\ar@{->>}[d]^{f'} \ar[u]|{s_0^{f'}} \\
Y \ar@{->>}[r]_{y} & Y'. }
$$
 The previous proposition shows that the non-dotted upper part of this diagram is a regular pushout, and consequently that the
 factorization $\psi \colon R[f] \to \overline{R[f']}$, where
 $\overline{R[f']}$ is the pullback of $(d_0^{f'},s_0^{f'})$ along $x$, is a regular epimorphism. The same
 method as in the previous proposition allows us to show that the factorization $\phi \colon X \to \bar X'$,  where $\bar X'$ is
 the pullback of $f'$ along $y$, is a regular epimorphism, and consequently that the square \eqref{diagram comm square of reg epis} is a
 regular pushout.
 \endproof

\section{The denormalized $3 \times 3$ lemma} \label{section denormalized 3x3 lemma}

\subsection{Preliminary observations}
A diagram:
$$\xymatrix{
    G  \ar@<-4pt>[r]_{d_0} \ar@<4pt>[r]^{d_1} & X \ar[l]|{s_0} \ar[r]^{f} & Y, }
$$
where the left-hand side part is a reflexive graph and $f$
coequalizes the pair $(d_0, d_1)$, is said to be \emph{left exact}
when $G$ is the kernel equivalence relation $R[f]$ of $f$, while
it is \emph{right exact} when $f$ is the coequalizer of the pair
$(d_0, d_1)$. It is said to be \emph{exact} when it is both
left and right exact. \\

A commutative diagram:
\begin{equation} \label{weakly denormalized 3x3 diagram}
\xymatrix@=25pt{
      R[\phi] \ar@<-1,ex>[d]_{d_0^{\phi}} \ar@<+1,ex>[d]^{d_1^{\phi}} \ar@<-1,ex>[r]_{R(d_1^x)}\ar@<+1,ex>[r]^{R(d_0^x)}
           & R[f] \ar@<-1,ex>[d]_{d_0^f}\ar@<+1,ex>[d]^{d_1^f} \ar[l] \ar[r]^{R(x)} & R[f'] \ar@<-1,ex>[d]_{d_0^{f'}}\ar@<+1,ex>[d]^{d_1^{f'}}\\
      R[x] \ar@{->>}[d]_{\phi}\ar[u] \ar@<-1,ex>[r]_{d_0^x}\ar@<+1,ex>[r]^{d_1^x} & X \ar@{->>}[d]^{f}\ar[u] \ar[l] \ar@{->>}[r]^{x} & X'
      \ar@{->>}[d]^{f'}\ar[u]\\
       W \ar@<-1,ex>[r]_{y_0} \ar@<+1,ex>[r]^{y_1}  & Y  \ar[l] \ar@{->>}[r]_{y} & Y'
                     }
\end{equation}
is said to be a \emph{weakly $3 \times 3$ diagram} when the middle
row and the three columns are exact, and a \emph{$3 \times 3$
diagram} when all the rows and columns are exact. In a weakly $3
\times 3$ diagram, the pair $(R(d_0^x),R(d_1^x))$ is necessarily
jointly monomorphic, and the map $y$ is necessarily an extremal
epimorphism. Since $\phi$ is an epimorphism, the lower row is
right exact if and only if the lower right-hand side square is a
pushout. Now suppose $\EE$ is regular.

\begin{prop}[\cite{Bourn denormalized 3x3}]\label{p2}
Let $\EE$ be a regular category. In a weakly $3 \times 3$ diagram,
the upper row is left exact if and only if the pair $(y_0,y_1)$ is
jointly monomorphic.
\end{prop}

We recall from \cite{Bourn denormalized 3x3} the ``denormalized $3
\times 3$ lemma'' for regular Mal'tsev categories:

\begin{prop}\label{p111}
Let $\EE$ be a regular Mal'tsev category. Given any weakly $3
\times 3$ diagram \eqref{weakly denormalized 3x3 diagram}, the
following conditions are equivalent:
\begin{itemize}
\item[(i)] the upper row is exact; \item[(ii)] the lower row is
exact; \item[(iii)] \eqref{weakly denormalized 3x3 diagram} is a
$3 \times 3$ diagram.
\end{itemize}
\end{prop}

In any category $\EE$, we shall say that a weakly $3 \times 3$
diagram \eqref{weakly denormalized 3x3 diagram} satisfies the
denormalized $3 \times 3$ lemma if the three previous conditions
are equivalent for the diagram.

\begin{prop}\label{p0}
Let $\EE$ be a regular category. A weakly $3 \times 3$ diagram
\eqref{weakly denormalized 3x3 diagram}, such that the lower
right-hand side square is a regular pushout, satisfies the  $3
\times 3$ lemma.
\end{prop}

\proof In such a diagram, the map $R(x)$ and $R(f)$ are
necessarily regular epimorphisms. Since the lower right-hand side
square is a pushout, the lower row is right exact. Suppose now
that the lower row is left exact; then the upper row is left exact
by Proposition \ref{p2}. Since $R(x)$ is a regular epimorphism,
the upper row is exact. Conversely, when the upper row is exact,
the pair $(y_0,y_1)$ is jointly monomorphic, and the factorization
$t \colon W \to R[y]$ is a monomorphism such that $t \,\phi =
R(f)$. Since $R(f)$ is a regular epimorphism, $t$ is a regular
epimorphism as well, and consequently an isomorphism. Accordingly
the lower row is left exact and, since $g$ is a regular
epimorphism, right exact.
\endproof

\subsection{The regular $\Sigma$-Mal'tsev context}

We shall now investigate what is remaining of the denormalized
$3\times 3$ lemma in a regular $\Sigma$-Mal'tsev category $\EE$.

\begin{theo}\label{th1}
Let $\EE$ be a regular $\Sigma$-Mal'tsev category. Consider a
weakly $3 \times 3$ diagram \eqref{weakly denormalized 3x3
diagram}:
\begin{enumerate}
\item if the morphism $x$ is $\Sigma$-special and the lower row is
exact, then the upper row is exact; \item if the morphism $f$ is
$\Sigma$-special and the upper row is exact, then the lower row is
exact; \item if both $f$ and $x$ are $\Sigma$-special, then the
diagram \eqref{weakly denormalized 3x3 diagram} satisfies the
denormalized $3\times 3$ lemma.
\end{enumerate}
\end{theo}

\proof
\begin{enumerate}
\item If the lower row is exact, we have $\phi=R(f)$. According to
Proposition \ref{p1}, the lower right-hand side square is a
regular pushout, since $x$ is $\Sigma$-special. Then the upper row
is exact by Proposition \ref{p0}. \item When the upper row is
exact, the map $R(x)$ is a regular epimorphim. When $f$ is
$\Sigma$-special, the lower right-hand side square is a regular
pushout. Again by Proposition \ref{p0}, the lower row is exact.
\item It is an immediate consequence of (1) and (2).
\end{enumerate}
\endproof

We shall now investigate the conditions under which the rows are
$\Sigma$-special. For that, let us recall from \cite{Bourn sigma
Maltsev} the following definition:

\begin{defi}
Let $\EE$ be a regular category, and $\Sigma$ a class of points.
This class is said to be \emph{$2$-regular} whenever, given any
regular epimorphism of points as on the right-hand side of the
following diagram:
$$\xymatrix{
R[x]   \ar@<-2pt>[d]_{R(f)}  \ar@<4pt>[r]^{d_0^x}
\ar@<-4pt>[r]_{d_1^x} &
X \ar[l] \ar[d]_{f} \ar@{->>}[r]^{x} & X' \ar@<-2pt>[d]_{f'}  \\
R[y]  \ar@<-2pt>[u]_{R(s)}  \ar@<4pt>[r]^{d_0^y}
\ar@<-4pt>[r]_{d_1^y} & Y \ar@{->>}[r]_{y} \ar[l]
\ar@<-4pt>[u]_{s} & Y' \ar@<-2pt>[u]_{s'} }
$$
the point $(f',s')$ is in $\Sigma$ as soon as both $(f,s)$ and
$(R(f),R(s))$ belong to $\Sigma$.
\end{defi}

\begin{prop}\label{examp}
In the category $Mon$ of monoids, the class $\Sigma$ of Schreier
points is $2$-regular (see \cite{Bourn Mal'tsev reflection} where
it is asserted, but not proved). More generally, given any
finitely complete regular category $\EE$, the class of internal
Schreier points in the category $Mon\EE$ of internal monoids is
$2$-regular.
\end{prop}

\proof Consider a horizontal morphism of split epimorphisms, as in
the lower right-hand side part of the following diagram:
$$\xymatrix{
K[R(f)] \ar@<-2pt>@{ >->}[d]_{k_{R(f)}} \ar@<4pt>[r]^{K(d_0^x)}
\ar@<-4pt>[r]_{K(d_1^x)} & K[f] \ar[l] \ar@<-2pt>@{
>->}[d]_{k_{f}}
\ar@{->>}[r]^{K(x)} & K[f'] \ar@<-2pt>@{ >->}[d]_{k_{f'}} \\
R[x] \ar@<-2pt>@{.>}[u]_{q_{R}} \ar@<-2pt>[d]_{R(f)}
\ar@<4pt>[r]^{d_0^x} \ar@<-4pt>[r]_{d_1^x} & X
\ar@<-2pt>@{.>}[u]_{q_{f}} \ar[l] \ar@<-2pt>[d]_{f} \ar@{->>}[r]^{x} & X' \ar@<-2pt>[d]_{f'} \ar@<-2pt>@{.>}@/_/[u]_{q_{f'}} \\
R[y] \ar@<-2pt>[u]_{R(s)} \ar@<4pt>[r]^{d_0^y}
\ar@<-4pt>[r]_{d_1^y} & Y \ar@{->>}[r]_{y} \ar[l]
\ar@<-2pt>[u]_{s} & Y'. \ar@<-2pt>[u]_{s'} }
$$
Complete the diagram with the horizontal kernel pairs and the
vertical kernels. Commutation of limits makes the upper row left
exact. Since $y$ and $x$ are regular epimorphims, and since
$(f,s)$ is a Schreier point, then the lower right-hand side square
is a regular pushout by Proposition \ref{p10}. Accordingly the map
$K(x)$ is a regular epimorphism and the upper row is right exact
as well. The map $x$ and $K(x)$ are still regular epimorphisms
(=surjective maps) in the category $Set$ of sets. Accordingly the
Schreier retractions $q_f$ and $q_R$ produce the desired
retraction $q_{f'}$. As for the internal case, it is easy to check
that $Mon\EE$ is regular when $\EE$ is. Then the proof is just an
internal version of the one we did for $Mon$.
\endproof

The class of acupuncturing points in the category $Qnd$ of
quandles is $2$-regular, see \cite{Bourn quandles}.

\begin{prop}\label{p3}
Let $\Sigma$ be a point-congruous and $2$-regular class of points
in a regular category $\EE$. Then the full subcategory $\Sigma
l/Y$ of the slice category $\EE/Y$ is a regular category.
Moreover, in the diagram \eqref{diagram proof p3} below, if $f$
factors as $f = g \, h$ with $h$ a regular epimorphism, the map
$g$ is $\Sigma$-special if and only if $f$ and $\phi$ are. In
other words, a map $h$ in $\EE/Y$ is a regular epimorphism in
$\Sigma l/Y$ if and only if it is a regular epimorphism in $\EE$
such that the kernel equivalence relation $R[f]$ belongs to
$\Sigma l/Y$. Accordingly, when $\EE$ is Barr-exact, the category
$\Sigma l/Y$ is Barr-exact as well.
\end{prop}

\proof Let us first show that $\Sigma l/Y$ admits quotients of
effective equivalence relations. Consider a morphism $h$ in
$\Sigma l/Y$ as in the lower right-hand side square in the diagram
below (which means that both $f$ and $f'$ are $\Sigma$-special):
\begin{equation} \label{diagram proof p3}
\xymatrix@=25pt{
      R[\phi] \ar@<-1,ex>[d]_{d_0^{\phi}} \ar@<+1,ex>[d]^{d_1^{\phi}} \ar@<-1,ex>[r]_{R(d_1^h)}\ar@<+1,ex>[r]^{R(d_0^h)}
           & R[f] \ar@<-1,ex>[d]_{d_0^f}\ar@<+1,ex>[d]^{d_1^f} \ar[l] \ar[r]^{R(h)} & R[f'] \ar@<-1,ex>[d]_{d_0^{f'}}\ar@<+1,ex>[d]^{d_1^{f'}}\\
      R[h] \ar[d]_{\phi}\ar[u] \ar@<-1,ex>[r]_{d_0^h}\ar@<+1,ex>[r]^{d_1^h} & X \ar[d]^{f}\ar[u] \ar[l] \ar[r]^{h} & X' \ar[d]^{f'}\ar[u]\\
       Y \ar@{=}[r]  & Y   \ar@{=}[r] & Y.
                     }
\end{equation}
By commutation of limits, the upper row is left exact, and so the
left-hand side of the two-level upper part of the diagram produces
the kernel equivalence relation of the right-hand side morphism
$(h,R(h))$ in the category $Eq\EE$ of equivalence relations in
$\EE$. Since $\Sigma$ is point-congruous and both $R[f]$ and
$R[f']$ are $\Sigma$-equivalence relations, so is $R[\phi]$, and
hence $\phi$ is $\Sigma$-special. Now let $h = m \, x$ be the
decomposition of $h$ with $m$ a monomorphism and $x$ a regular
epimorphism, so that we have $R[h]=R[x]$ and $R[\phi]=R[R(x)]$. In
this way, we get the following regular epimorphism $(x,R(x))$
between equivalence relations on the right-hand side:
  $$ \xymatrix@=25pt{
        R[\phi]=R[R(x)] \ar@<-1,ex>[d]_{d_0^{\phi}} \ar@<+1,ex>[d]^{d_1^{\phi}} \ar@<-1,ex>[r]_-{R(d_1^h)}\ar@<+1,ex>[r]^-{R(d_0^h)}
             & R[f] \ar@<-1,ex>[d]_{d_0^f}\ar@<+1,ex>[d]^{d_1^f} \ar[l] \ar@{->>}[r]^{R(x)} & R[f' \, m] \ar@<-1,ex>[d]_{d_0^{f'}}\ar@<+1,ex>[d]^{d_1^{f'}}\\
        R[h]=R[x] \ar[u] \ar@<-1,ex>[r]_-{d_0^h}\ar@<+1,ex>[r]^-{d_1^h} & X \ar[u] \ar[l] \ar@{->>}[r]_{x} & X'', \ar[u]                       }
    $$
because the morphism $R(x)$ is a regular epimorphism: indeed, by
$(f' \, m) \, x = f$, $x$ is a regular epimorphism in the regular
category $\EE/Y$. Then so is the product $x\times_Yx=R(x)$ in this
same category. According to the $2$-regularity of $\Sigma$, the
right-hand side equivalence relation is a $\Sigma$-equivalence
relation, and $x \colon f \to f' \, m$ is a regular epimorphism in
$\Sigma l/Y$. Since $\EE$ is regular, it is clear that these
regular epimorphisms in $\Sigma l/Y$ are stable under pullbacks.
Suppose now that $\EE$ is Barr-exact and $R$ is an equivalence
relation on the object $f$ in $\Sigma l/Y$, which means that the
map $\phi=f \, d_0^R=f \, d_1^R$ is $\Sigma$-special. Since $\EE$
is Barr-exact, there is some regular epimorphism $h$ in $\EE$ such
that $R=R[h]$, whence the factorization $f'$ as above, with $f'$
$\Sigma$-special, according to the first  part of this proof.
Hence the equivalence relation $R$ is effective in $\Sigma l/Y$.
\endproof

From now on, saying that an exact fork is \emph{$\Sigma$-exact}
will mean that its regular epimorphic part is $\Sigma$-special.

\begin{theo}\label{th2}
Let $\EE$ be a regular $\Sigma$-Mal'tsev category. Consider a
weakly $3 \times 3$ diagram \eqref{weakly denormalized 3x3
diagram}. Then:
\begin{enumerate}
\item if $\Sigma$ is point-congruous, the map $x$ is
$\Sigma$-special and the lower row is $\Sigma$-exact, then the
upper row is $\Sigma$-exact; \item if the class $\Sigma$ is
$2$-regular, the maps $f$, $x$ are in $\Sigma$ and the upper row
is $\Sigma$-exact, then the lower row is $\Sigma$-exact.
\end{enumerate}
\end{theo}

\proof
\begin{enumerate}
\item We know, by Theorem \ref{th1}, that the upper row is exact.
Moreover, since $\Sigma$ is point-congruous, the category $\Sigma
l\EE$ is stable under finite limits in $\EE^2$. Accordingly, since
both $x$ and $y$ are $\Sigma$-special, so is $R(x)$. \item The map
$f$ being $\Sigma$-special and the upper row being exact, the
lower row is exact by Theorem \ref{th1}. Since $x$ and $R(x)$ are
$\Sigma$-special, $R[x]$ and $R[\phi] = R[R[x]]$ are
$\Sigma$-equivalence relations. Since $\Sigma$ is $2$-regular, $W
= R[y]$ is a $\Sigma$-equivalence relation, too. Consequently, the
map $y$ is $\Sigma$-special.
\end{enumerate}
\endproof

There is a last situation dealing with the denormalized $3 \times
3$ lemma. We noticed that $\Sigma l/Y$ is a finitely complete
category when $\Sigma$ is point-congruous, and a Mal'tsev category
when $\EE$ is a $\Sigma$-Mal'tsev category. By Proposition
\ref{p3}, we know moreover that $\Sigma l/Y$ is regular when
$\Sigma$ is $2$-regular.

\begin{prop}\label{p18}
Let $\EE$ be a regular $\Sigma$-Mal'tsev category such that the
class $\Sigma$ is point-congruous and $2$-regular. Any weakly
$3\times 3$ diagram \eqref{weakly denormalized 3x3 diagram} in
$\EE$ such that $y$, $f'$ and $f' \, x = y \, f$ are
$\Sigma$-special satisfies the denormalized $3\times 3$ lemma.
\end{prop}

\proof Let us think of diagram \eqref{weakly denormalized 3x3
diagram} as a diagram in the regular category $\EE/Y'$. In this
diagram, from any object there is a unique map to $Y'$, so that
any object in $\EE/Y'$ can be identified with its domain. Here any
object lies in $\Sigma l/Y'$ except $R[\phi]$ and $W$.  When the
lower row is exact, $W$ is in $\Sigma l/Y'$, and so is $R[\phi]$.
When the upper row is exact, $R[\phi]$ belongs to $\Sigma l/Y'$,
and according to Proposition \ref{p3} so does $W$. So, under any
of the conditions of the denormalized $3\times 3$ lemma, the whole
diagram lies in the regular Mal'tsev category $\Sigma l/Y'$, and
the denormalized $3\times 3$ lemma holds by Proposition
\ref{p111}.
\endproof

\section{$\Sigma$-protomodular categories} \label{section Sigma-protomodular categories}

\subsection{Preliminary observations}

We recall that a category $\EE$ is said to be \emph{protomodular}
\cite{Bourn protomod} when any change-of-base functor with respect
to the fibration of points is conservative. In the pointed
context, this condition implies that the category $\EE$ shares
with the category $Gp$ of groups the following well-known four
properties:
\begin{itemize}
\item[(i)] a morphism $f$ is a monomorphism if and only if its
kernel $K[f]$ is trivial; equivalently, pulling back reflects
monomorphisms;

\item[(ii)] any regular epimorphism is the cokernel of its kernel;
in other words, any regular epimorphism produces an exact
sequence;

\item[(iii)] there is a very specific class of monomorphisms $m
 \colon U \into X$, the normal ones, namely those such that there
exists a (necessarily unique) equivalence relation $R$ on $X$ such
that $m^{-1}(R)$ is the indiscrete equivalence relation on $X$ and
any commutative square in the following induced diagram is a
pullback:
$$\xymatrix{
U \times U \ar@<-6pt>[d]_{d_0^U}  \ar@{ >->}[r]^-{\check m}
\ar@<6pt>[d]^{d_1^U} & R \ar@<-6pt>[d]_{d_0^R}
\ar@<6pt>[d]^{d_1^R} \\
U \ar@{ >->}[r]_{m} \ar[u]|{s_0^U} & X; \ar[u]|{s_0^R} }
$$

\item[(iv)] any reflexive relation is an equivalence relation,
i.e. the category $\EE$ is a Mal'tsev one.
\end{itemize}

We are now interested in seeing what remains of the properties
above if we consider categories that are protomodular relatively
to a pullback stable class $\Sigma$ of points. To do that, let us
first recall the following definition, see \cite{Bourn monad
intgpd} and \cite{MMS semidir normal cats}:

\begin{defi}
A point $(f,s)$ in $\EE$ is called \emph{strong} whenever, given
any pullback:
$$
\xymatrix{
\bar X \pullback \ar[r]^{x} \ar@<-2pt>[d]_{\bar f} & X \ar@<-2pt>[d]_{f} \\
\bar Y \ar[r]_{y}  \ar@<-2pt>[u]_{\bar s} & Y, \ar@<-2pt>[u]_{s} }
$$
the pair $(x,s)$ is jointly extremally epimorphic.
\end{defi}

and also the following one, see \cite{Bourn Mal'tsev reflection},
and \cite{BMMS S-protomodular}\label{sproto} in the pointed case:

\begin{defi}\label{def1}
Let $\EE$ be a category endowed with a pullback stable class
$\Sigma$ of points. $\EE$ is said to be
\emph{$\Sigma$-protomodular} when every point in $\Sigma$ is
strong.
\end{defi}

A category $\EE$ is protomodular if and only if every point is
strong. Let us now review the four previous properties w.r.t this
concept of partial protomodularity. As for (iv), we have the
following:

\begin{theo}[\cite{Bourn Mal'tsev reflection}] \label{main}
Let $\EE$ be a category endowed with a pullback stable class
$\Sigma$ of points. Then:
\begin{enumerate}
\item when $\EE$ is $\Sigma$-protomodular, it is a
$\Sigma$-Mal'tsev category;

\item when, in addition, $\Sigma$ is point-congruous, any
change-of-base functor with respect to the subfibration
$\P_{\EE}^{\Sigma}$ of the fibration of points is conservative.
\end{enumerate}
\end{theo}

As for (i) we get:

\begin{prop}
Let $\EE$ be a $\Sigma$-protomodular category. Then:
\begin{enumerate}
\item the split epimorphic part of a point in $\Sigma$ is an
isomorphism if and only if any of its pullbacks is an isomorphism;

\item pulling back $\Sigma$-special morphisms reflects
monomorphisms.
\end{enumerate}
In particular, when $\EE$ is pointed, more classically we get:
\begin{enumerate}
\item the split epimorphic part of a point in $\Sigma$ is an
isomorphism if and only if its kernel is trivial;

\item a $\Sigma$-special morphism is a monomorphism if and only if
its kernel is trivial.
\end{enumerate}
\end{prop}

\proof
\begin{enumerate}
\item Consider the following pullback, with $(f,s)$ in $\Sigma$,
and suppose $f'$ is an isomorphism:
  $$\xymatrix{
      X' \pullback \ar[r]^{\bar g} \ar@<-2pt>[d]_{f'} & X  \ar@<-2pt>[d]_f\\
      Y'   \ar[r]^g \ar@<-2pt>[u]_{s'} & Y.  \ar@<-2pt>[u]_s }
    $$
We get $s' \, f' = 1_{X'}$. We can check that $s \, f = 1_X$ by
composing with the jointly extremally epimorphic pair $(s,\bar
g)$. The equality $s \, f \, s = 1_X$ is straightforward.
Moreover, $s \, f \, \bar g = \bar g \, s' \, f' = \bar g = 1_X \,
\bar g$.

\item Consider the following diagram, where $f$ is
$\Sigma$-special and the lower square is a pullback:
$$\xymatrix{
R[f'] \ar@<-6pt>[d]_{d_0^{f'}} \ar[r]^{R(\bar g)} \ar@<6pt>[d]^{d_1^{f'}} & R[f] \ar@<-6pt>[d]_{d_0^f} \ar@<6pt>[d]^{d_1^f} \\
X' \pullback \ar[r]_{\bar g} \ar[d]_{f'} \ar[u]|{s_0^{f'}} & X
\ar[d]^f
\ar[u]|{s_0^f} \\
Y' \ar[r]_g & Y. }
$$
Then the two upper commutative squares are pullbacks, and the pair
$(d_0^f,s_0^f)$ is in $\Sigma$. Suppose moreover that $f'$ is a
monomorphism; then $d_0^{f'}$ is an isomorphism. According to (1),
the map $d_0^f$ is itself an isomorphism, and $f$ is a
monomorphism.
\end{enumerate}
\endproof

As for (ii), we get:
\begin{prop}\label{p5}
Let $\EE$ be a $\Sigma$-protomodular category. Then:
\begin{enumerate}
\item given any pullback of the point $(f,s)$ in $\Sigma$:
$$\xymatrix{
X' \pullback \ar[r]^{\bar g} \ar@<-2pt>[d]_{f'} & X  \ar@<-2pt>[d]_f \\
Y' \ar[r]_g \ar@<-2pt>[u]_{s'} & Y, \ar@<-2pt>[u]_s }
      $$
   the downward square is a pushout;

\item consider any pullback of the form:
$$\xymatrix{
X' \pullback \ar[r]^{\bar g} \ar[d]_{f'} & X  \ar[d]^f \\
Y'  \ar[r]_g & Y. }
$$
If there exists a factorization $f \, h = g$ and $f$ is a
$\Sigma$-special regular epimorphism, the square is a pushout as
well.
\end{enumerate}
In particular, when $\EE$ is pointed, more classically we get:
\begin{enumerate}
\item any split epimorphism in $\Sigma$ is the cokernel of is
kernel;

\item any $\Sigma$-special regular epimorphism is the cokernel of
its kernel.
\end{enumerate}
\end{prop}

\proof
\begin{enumerate}
\item Suppose we have a pair of morphisms $\alpha \colon X \to Z$,
$\gamma \colon Y' \to Z$ such that $\alpha \, \bar g = \gamma \,
f'$. We have that $\alpha \, s \, f = \alpha$: this can be checked
by composing with the jointly extremally epimorphic pair $(s, \bar
g)$. Since $f$ is an epimorphism, this factorization is unique. It
remains to check that $\alpha \, s \, g = \beta$; this can be done
by composing with the epimorphism $f'$.

\item Consider the same pullback without splittings, with $f$
$\Sigma$-special and $f \, h = g$. Then the map $h$ produces a
splitting $(1, h) \colon Y' \to X'$ of $f$', which is then an
epimorphism. Now complete the diagram with the kernel equivalence
relations:
$$\xymatrix{
R[f'] \ar@<-6pt>[d]_{d_0^{f'}} \ar[r]^{R(\bar g)} \ar@<6pt>[d]^{d_1^{f'}} & R[f] \ar@<-6pt>[d]_{d_0^f} \ar@<6pt>[d]^{d_1^f} \\
X' \pullback \ar[d]_{f'} \ar[r]_{\bar g} \ar[u]|{s_0^{f'}} & X
\ar[d]^f \ar[u]|{s_0^f} \\
Y'  \ar[r]_g & Y.  }
$$
Any commutative square in the upper part is a pullback; moreover,
the pair $(d_0^f, s_0^f)$ is in $\Sigma$, because $f$ is
$\Sigma$-special. Suppose we have a pair of morphisms $\alpha
\colon X \to Z$, $\gamma \colon Y' \to Z$ such that $\alpha \,
\bar g = \gamma \, f'$. Let us show that $\alpha$ coequalizes the
pair $(d_0^f,d_1^f)$. This can be checked by composing with the
jointly extremally epimorphic pair $(s_0^f,R(\bar g))$. We have
\[ \alpha \, d_0^f \, s_0^f = \alpha = \alpha \, d_1^f \, s_0^f, \] and
\[ \alpha \, d_0^f \, R(\bar{g}) = \alpha \, \bar{g} \, d_0^{f'} =
\gamma \, f' \, d_0^{f'} = \gamma \, f' \, d_1^{f'} = \alpha \,
\bar{g} \, d_1^{f'} = \alpha \, d_1^f \, R(\bar{g}). \] Since $f$
is a regular epimorphism, there is a unique factorization
$\bar{\alpha} \colon Y \to Z$ such that $f \, \bar{\alpha} =
\alpha$. The equality $\bar{\alpha} \, g = \gamma$ can be checked
by composing with the epimorphism $f'$.
\end{enumerate}
\endproof

In the regular context we can add the following precision:

\begin{prop}\label{p4}
Let $\EE$ be a pointed regular $\Sigma$-protomodular category and
$f$ a $\Sigma$-special regular epimorphism. Let $g \colon X \to
Y'$ be a $\Sigma$-special morphism such that $k_g = k_f$,
 where $k_g$ and $k_f$ are the kernels of $g$ and $f$,
respectively. Then the unique factorization $t$ such that $g=t
\,f$ is a monomorphism.
\end{prop}

\proof We have $0=g \, k_g=g \, k_f$. Moreover $f$, being
$\Sigma$-special, is the cokernel of $k_f$ (by Proposition
\ref{p5}). Hence we get the factorizarion $t$, and an inclusion $i
\colon R[f] \rightarrowtail R[g]$ between two $\Sigma$-equivalence
relations. Now consider the following diagram:
$$\xymatrix{
             & {R[f]\;} \ar@{>->}[rrd]^i \ar@<-2pt>[ddr]_>>>>>>{d_0^f}\\
           {K[f]=K[g]\;}  \ar@{>->}[rrr]^<<<<<{(0,k_g)}\ar@{>->}[ru]^{(0,k_f)} \ar@<-2pt>[d]_{f'} &&& R[g]  \ar@<-2pt>[d]_{d_0^g}\\
           {1\;}   \ar@{>->}[rr] \ar@<-2pt>[u] && Y   \ar@<-2pt>[uul]_<<<<<<{s_0^f} \ar@{=}[r] & Y. \ar@<-2pt>[u]_{s_0^g } }
$$
The rectangle is a pullback and $(d_0^g,s_0^g)$ in $\Sigma$, so
the pair $(s_0^g,(0,k_g))$ is jointly extremally epimorphic. Since
the pair $(s_0^f,(0,k_f))$ factors through the monomorphism $i$,
$i$ is an isomorphism. Hence we get $R[f] \simeq R[g]$, which
implies, in the regular category $\EE$, that $t$ is a
monomorphism.
\endproof

As for (iii), we recall from \cite{Bourn Mal'tsev reflection} the
following:
\begin{prop}\label{normal}
Let $\EE$ be a $\Sigma$-Mal'tsev category. When $m \colon U
\rightarrowtail X$ is a monomorphism which is normal to a
$\Sigma$-equivalence relation $S$ on $X$, the object $U$ is
$\Sigma$-special. When $\EE$ is $\Sigma$-protomodular, a
monomorphism $m$ is normal to at most one $\Sigma$-equivalence
relation.
\end{prop}

As a first step toward the $3\times 3$ lemma, we get:
\begin{prop}
Let $\EE$ be a pointed $\Sigma$-protomodular category. Consider
any morphism between points, with codomain $(f,s)$ in $\Sigma$:
\begin{equation} \label{diagram morphism between points}
\xymatrix{
X' \ar[r]^{\bar g} \ar@<-2pt>[d]_{f'} & X  \ar@<-2pt>[d]_f \\
Y' \ar[r]_g \ar@<-2pt>[u]_{s'} & Y . \ar@<-2pt>[u]_s }
\end{equation}
Suppose moreover that the factorization $\phi \colon X' \to \check
X$, where $\check X$ is the pullback of $(f,s)$ along $g$, is a
monomorphism. Then the square \eqref{diagram morphism between
points} is a pullback if and only if the map $K(\bar g) \colon
K[f'] \to K[f]$ is an isomorphism. It is clear that the condition
on $\phi$ holds as soon as $\bar g$ is a monomorphism (the
splittings making then $g$ itself a monomorphism).
\end{prop}

\proof Clearly {if \eqref{diagram morphism between points} is a
pullback, then $K(\bar g)$ is an isomorphism. Conversely,}
consider the following diagram with exact columns:
$$\xymatrix{  K[f'] \ar[r]^{K(\bar g)}\ar@{ >->}[d]_{k_{f'}}  & K[\check f] \ar@{=}[r]\ar@{ >->}[d]_{(0,k_{f})} & K[f] \ar@{ >->}[d]_{k_{f}} \\
     X' \ar@<-2pt>[d]_{f'} \ar@{>->}[r]^{\phi}  &  \check X \pullback \ar[r]^{\check  g} \ar@<-2pt>[d]_{\check f} & X  \ar@<-2pt>[d]_f\\
     Y'\ar@<-2pt>[u]_{s'} \ar@{=}[r] &  Y'   \ar[r]^g \ar@<-2pt>[u]_{\check s} & Y.  \ar@<-2pt>[u]_s }
$$
Since $\EE$ is $\Sigma$-protomodular and $(f,s)$ in $\Sigma$,
$(\check{f}, \check{s})$ is in $\Sigma$, too, and the pair
$(\check s, (0,k_{f}))$ is jointly extremally epimorphic. Since
$s'$ and $k_{f'} \, K(\bar g)^{-1}$ factor through the
monomorphism $\phi$, $\phi$ is an isomorphism, and the square
\eqref{diagram morphism between points} is a pullback.
\endproof

\begin{coro}\label{cor1}
Let $\EE$ be a pointed regular $\Sigma$-protomodular category.
\begin{enumerate}
\item Let $f'$ be a regular epimorphism such that $f'=f \, m$ with
$m$ a monomorphism and $f$ $\Sigma$-special. If $k(m) \colon K[f']
\to K[f]$ is an isomorphism, then $m$ is an isomorphism, too.

\item Consider a commutative diagram:
\begin{equation} \label{diagram cor 4.10}
\xymatrix{  {K[f']\;} \ar@{->>}[d]_{K(\bar g)} \ar@{>->}[r]^{k_{f'}}  & X'  \ar@{->>}[r]^{f'} \ar[d]_{\bar g} & Y'  \ar@{->>}[d]^g\\
       {K[f]\;} \ar@{>->}[r]_{k_f} &  X  \ar@{->>}[r]_f  & Y }
\end{equation}
where $f$ is a $\Sigma$-special regular epimorphism and $f'$ is a
regular epimorphism. Then $\bar g$ is a regular epimorphism and
the right-hand side square is a regular pushout as soon as both
$g$ and $K(\bar g)$ are regular epimorphisms. Accordingly the map
$K(f') \colon K[\bar g] \to K[g]$ is a regular epimorphism as
well.
\end{enumerate}
\end{coro}

\proof
\begin{enumerate}
\item Consider the following diagram:
$$\xymatrix{
K[f']\; \ar@{ >->}[d]_{K(m)}^{\simeq} \ar@{>->}[r]^{(0, k_{f'})} &
R[f'] \ar@{ >->}[d]_{R(m)} \ar@<4pt>[r]^-{d_0^{f'}}
\ar@<-4pt>[r]_-{d_1^{f'}} & X' \ar[l]|-{s_0^{f'}} \ar@{ >->}[d]_{m} \ar@{->>}[r]^{f'} & Y \ar@{=}[d]  \\
K[f]\; \ar@{>->}[r]_{(0, k_{f})} & R[f] \ar@<4pt>[r]^{d_0^f}
\ar@<-4pt>[r]_{d_1^f} & X \ar@{->>}[r]_{f} \ar[l]|{s_0^f} & Y. }
$$
We can apply the previous proposition to the central square,
because $K(m)$ is an isomorphism. The point $(d_0^f, s_0^f)$ is in
$\Sigma$, since $f$ is $\Sigma$-special. Consequently the central
square indexed by $0$ a pullback, and this central part of the
diagram becomes a discrete fibration between equivalence
relations. According to the Barr-Kock Theorem valid in any regular
category (see, for example, \cite{BG categorical foundations}),
the right-hand side square is a pullback as well, since $f$ and
$f'$ are regular epimorphisms. Accordingly $m$ is an isomorphism.

\item Consider the diagram \eqref{diagram cor 4.10}. Let $\bar g =
m \, \check g$ be the decomposition of $\bar g$ with $m \colon U
\into X$ a monomorphism and $\check g$ a regular epimorphism. Then
$m \, f \colon U \to Y$ is a regular epimorphism, because $g$ and
$f'$ are. Since $K(\bar g)$ is a regular epimorphism, and hence a
strong epimorphism, there is a unique factorization $k \colon K[f]
\to U$ such that $m \, k = k_f$, and this $k$ is necessarily the
kernel of the regular epimorphism $f \, m$. Whence
$K(m)=1_{K[f]}$, so that, according to the first point, $m$ is an
isomorphism and $\bar g$ is a regular epimorphism.

Now take the pullback $\bar f$ of $f$ along $g$:
$$
\xymatrix{ \bar{X} \ar[r]^{\bar{f}} \ar[d]_{\check{g}} & Y'
\ar[d]^g \\
X \ar[r]_f & Y. }
$$
$\bar{f}$ is a $\Sigma$-special regular epimorphism, since $f$ is.
Let $\bar g = \check g \, \psi$ be the induced decomposition
through this pullback. Consider now the following commutative
diagram, where the lower row is a $\Sigma$-special exact sequence
and $f'$ is a regular epimorphism:
$$\xymatrix{  {K[f']\;} \ar@{->>}[d]_{K(\bar g)=K(\psi)} \ar@{>->}[r]^{k_{f'}}  & X'  \ar@{->>}[r]^{f'} \ar[d]_{\psi} & Y  \ar@{=}[d]\\
          {K[f]\;} \ar@{>->}[r]_{(k_f,0)} &  \bar X  \ar@{->>}[r]_{\bar f}  & Y } $$
The map $\psi$ is a regular epimorphism, since $K(\psi)=K(\bar g)$
is a regular epimorphism. Accordingly, the square $g \, f' = f \,
\bar{g}$ is a regular pushout.
\end{enumerate}
\endproof

\section{Aspects of the $3\times 3$ lemma} \label{section aspects 3x3 lemma}

In a pointed category, an exact sequence is a sequence $\xymatrix{
K[f] \ar@{ >->}[r]^-{k_f} & X \ar@{->>}[r]^f & Y}$ where $k_f$ is
the kernel of $f$ and $f$ the cokernel of $k_f$. In a pointed
protomodular category, a regular epimorphism $f$ is the cokernel
of its kernel, and  so it gives rise to an exact sequence. In a
$\Sigma$-protomodular category, Proposition \ref{p5} shows that
any $\Sigma$-special regular epimorphism $f$ determines an exact
sequence, as well.

A commutative diagram:
\begin{equation} \label{weakly normalized 3x3 diagram}
\xymatrix@=25pt{
      {K[\phi]\;} \ar@{ >->}[d]_{k_{\phi}}  \ar@{>->}[r]^{K(k_x)} & K[f] \ar@{ >->}[d]^{k_f} \ar[r]^{K(x)} & K[f'] \ar@{ >->}[d]^{k_{f'}}\\
      {K[x]\;} \ar@{->>}[d]_{\phi} \ar@{>->}[r]_{k_x} & X \ar@{->>}[d]^{f} \ar@{->>}[r]_{x} & X' \ar@{->>}[d]^{f'}\\
      U \ar[r]_{u}  & Y   \ar@{->>}[r]_{y} & Y' }
\end{equation}
in a pointed category $\EE$ is said to be a \emph{weakly $3\times
3$ diagram} when the three columns and the middle row are exact
and a \emph{$3\times 3$ diagram} when all the columns and rows are
exact. In a pointed regular protomodular category the
\emph{$3\times 3$ lemma} holds:

\begin{prop}[\cite{Bourn 3x3}]
Given a weakly $3\times 3$ diagram \eqref{weakly normalized 3x3
diagram} in a pointed regular protomodular category $\EE$, the
following conditions are equivalent:
\begin{enumerate}
\item the upper row is exact;

\item the lower row is exact;

\item \eqref{weakly normalized 3x3 diagram} is $3\times 3$
diagram.
\end{enumerate}
\end{prop}

Given a weakly $3\times 3$ diagram \eqref{weakly normalized 3x3
diagram} in a pointed category, we say that it satisfies the
$3\times 3$ lemma when the three previous conditions are
equivalent for the diagram.

Suppose now that $\EE$ is a pointed regular $\Sigma$-protomodular
category. We shall be interested in the \emph{weakly
$\Sigma$-special $3\times 3$ diagrams}, namely those diagrams
whose three columns and middle row are $\Sigma$-exact, i.e. where
$f$, $f'$, $\phi$ and $x$ are $\Sigma$-special regular
epimorphisms, and consequently produce exact sequences. \\

The next proposition is a version for $\Sigma$-protomodular
categories of the so-called \emph{upper $3 \times 3$ lemma}:

\begin{prop}\label{p22}
Let $\EE$ be a pointed regular $\Sigma$-protomodular category and
let \eqref{weakly normalized 3x3 diagram} be a weakly
$\Sigma$-special $3\times 3$ diagram. Suppose moreover that
$\Sigma$ is point-congruous. When the lower row is $\Sigma$-exact,
so is the upper row.
\end{prop}

\proof By assumption, $y$ is a $\Sigma$-special regular
epimorphism, $u$ is the kernel of  $y$ and $\phi = K(f)$. This
last equality implies that the upper row is left exact by
commutation of limits, namely that $K(k_x)$ is the kernel of
$K(x)$.

According to Corollary \ref{cor1}, this implies that the lower
right-hand side square in \eqref{weakly normalized 3x3 diagram} is
a regular pushout, and hence $K(x)$ is a regular epimorphism. When
$\Sigma$ is point-congruous, the map $K(x)$ is in $\Sigma$,
because $x$ and $y$ are.
\endproof

The implication in the inverse direction (i.e. the so-called
\emph{lower $3 \times 3$ lemma}) needs much more information, and
it points out a totally asymmetric situation.

\begin{defi}
Let $\Sigma$ be a pullback stable class of points in a pointed
category $\EE$. We say that this class is \emph{equi-consistent}
when, given any split epimorphism of equivalence relations:
\begin{equation} \label{diagram equi-consistent}
\xymatrix@=25pt{
      R \ar@<-1,ex>[d]_{d_0^{R}} \ar@<+1,ex>[d]^{d_1^{R}} \ar@<-1,ex>[r]_{\bar g}
           & S \ar@<-1,ex>[d]_{d_0^S}\ar@<+1,ex>[d]^{d_1^S} \ar[l]_{\bar t}  \\
      X \ar[u] \ar@<-1,ex>[r]_{g} &  Y\ar[u] \ar[l]_t  }
\end{equation}
where $R$ a $\Sigma$-equivalence relation and the split
epimorphism $(g,t)$ is in $\Sigma$, $(\bar g,\bar t)$ is in
$\Sigma$ as soon as $(K_0(\bar g),K_0(\bar t))$ is in $\Sigma$,
where $K_0(\bar g)$ and $K_0(\bar t)$ are the restrictions of
$\bar{g}$ and $\bar{t}$ to the kernels of $d_0^R$ and $d_0^S$.
\end{defi}

When $\Sigma$ is point-congruous, the equivalence relation $S$ is
a $\Sigma$-one as well, since the morphism $(t,\bar t)$ is an
equalizer in $Pt\EE$. Hence the previous condition becomes a
characteristic one: namely, under the same assumptions of the
previous proposition, the point $(\bar g,\bar t)$ is in $\Sigma$
if and only if the point $(K_0(\bar g),K_0(\bar t))$ is in
$\Sigma$. Indeed, being $\Sigma$ a point-congruous class, the
kernel of a map in $\Sigma Pt\EE$ is still in $\Sigma Pt\EE$.

\begin{prop}\label{e-q}
The class of Schreier points in $Mon$ is equi-consistent. The same
holds for the class of internal Schreier points in $Mon\EE$.
\end{prop}

\proof Consider a commutative diagram like \eqref{diagram
equi-consistent} in $Mon$. Denote by $q$ the Schreier retraction
of the Schreier point $(g,t)$; thanks to Proposition $2.4$ in
\cite{BMMS SemForum}, we know that the map $q$ satisfies the
equalities $x=q(x) \cdot tg(x)$ for all $x\in X$, and $q(k \cdot
t(y))=k$ for all $(k,y) \in K[g] \times Y$. Proposition $2.1.5$ in
\cite{Schreier book} tells us that the following two equalities
also hold:
\begin{enumerate}
\item $q(x \cdot x')=q(x) \cdot q(tg(x) \cdot q(x')) \qquad
\text{for all} \, (x,x') \in X\times X$;

\item $q(t(y) \cdot k) \cdot t(y)=t(y) \cdot k \qquad \text{for
all} (k,y) \in K[g] \times Y$.
\end{enumerate}

Because of the uniqueness of $q$, it suffices to show that
$q(x)Rq(x')$ whenever $xRx'$. The equi-consistent assumption means
that, when $x=1$, we have $1Rq(x')$. Since $R$ is a Schreier
equivalence relation, there is a Schreier retraction $\chi \colon
R \to K[d_0^R]$ satisfying the equalities $1R\chi(xRx')$,
$x'=\chi(xRx') \cdot x$ and $\chi(zR(u \cdot z))=u$ for all $z\in
X$, whenever we have $1Ru$. Starting with $xRx'$, we get
$1R\chi(xRx')$ and $x'=\chi(xRx') \cdot x'$. Hence we get from (1)
that:
$$q(x')=q(\chi(x,x')) \cdot q(tg\chi(x,x') \cdot q(x)).$$
So, from $1Rq\chi(xRx')$, we get also $(q(tg\chi(x,x') \cdot
q(x)))Rq(x')$. Now, from (2) we have:
$$q(tg\chi(x,x') \cdot q(x)) \cdot tg\chi(x,x')=tg\chi(x,x') \cdot q(x).$$
So, from $1Rtg\chi(x,x')$, we get $$q(x)Rtg\chi(x,x') \cdot q(x)
\, \text{and} \, q(tg\chi(x,x') \cdot q(x))Rtg\chi(x,x') \cdot
q(x).$$ Whence the following situation:
$$\xymatrix{
          q(x)   \ar[r]^R \ar@{.>}[rd] & tg\chi(x,x') \cdot q(x)\\
         q(tg\chi(x,x') \cdot q(x))  \ar[r]_{R} \ar[ru]^<<<<<<{R}  & q(x'), }
$$
and $q(x)Rq(x')$ as desired. The second assertion is a
straightforward consequence of the Yoneda embedding.
\endproof

\begin{prop}\label{p6}
Let $\Sigma$ be a pullback stable, $2$-regular, equi-consistent
class of points in a pointed regular category $\EE$. Consider any
regular epimorphism in $Pt\EE$:
$$\xymatrix{
         X'  \ar@{->>}[r]^{\bar g} \ar@<-2pt>[d]_{f'} & X  \ar@<-2pt>[d]_f\\
         Y'   \ar@{->>}[r]_g \ar@<-2pt>[u]_{s'} & Y  \ar@<-2pt>[u]_s }
       $$
with domain $(f',s')$ in $\Sigma$ and the morphisms $g$ and $\bar
g$ $\Sigma$-special. Then $(f,s)$ is in $\Sigma$ as soon as the
restriction $(K(f'),K(s'))$ to the kernels of $g$ and $\bar{g}$ is
in $\Sigma$.
\end{prop}

\proof Consider the following diagram:
$$\xymatrix{
   {K[\bar g]\;} \ar@<-2pt>[d]_{K(f')} \ar@{>->}[r]^{(0, k_{\bar g})}  &  R[\bar g] \ar@<-2pt>[d]_{R(f')} \ar@<1ex>[r]^{d_0^{\bar g}}
   \ar@<-1ex>[r]_{d_1^{\bar g}}  & X'  \ar@{->>}[r]^{\bar g} \ar@<-2pt>[d]_{f'}\ar[l] & X  \ar@<-2pt>[d]_f\\
   {K[g]\;} \ar@<-2pt>[u]_{K(s')} \ar@{>->}[r]_{(0, k_{g})}  &  R[g] \ar@<-2pt>[u]_{R(s')}  \ar@<1ex>[r]^{d_0^{g}} \ar@<-1ex>[r]_{d_1^{g}}
    &   Y'   \ar@{->>}[r]_g \ar@<-2pt>[u]_{s'} \ar[l] & Y.  \ar@<-2pt>[u]_s }
$$
Since $R[g]$ and $R[\bar{g}]$ are $\Sigma$-equivalence relations
and the class $\Sigma$ is equi-consistent, the pair
$(R(f'),R(s'))$ is in $\Sigma$. Moreover, since $\Sigma$ is
$2$-regular, the pair $(f,s)$ is in $\Sigma$ as well.
\endproof

\begin{coro}\label{cor2}
Let $\Sigma$ be a pullback stable, point-congruous, $2$-regular,
equi-consistent class of points in a pointed regular category
$\EE$. Consider any regular pushout:
$$\xymatrix{
         X'  \ar@{->>}[r]^{x} \ar@{->>}[d]_{f'} & X  \ar@{->>}[d]^f\\
         Y'   \ar@{->>}[r]_y  & Y }
$$
such that the morphisms $f'$, $y$ and $x$ are $\Sigma$-special.
Then the map $f$ is $\Sigma$-special as soon as the restriction
$K(f') \colon K[x]\to K[y]$ is.
\end{coro}

\proof Consider the following diagram:
$$\xymatrix{
{R[K(f')]\;}
\ar@<-1,ex>[d]_{d_0^{K(f')}}\ar@<+1,ex>[d]^{d_1^{K(f')}}
\ar@{>->}[r]^-{R(k_{x})} & R[f']
\ar@<-1,ex>[d]_{d_0^{f'}}\ar@<+1,ex>[d]^{d_1^{f'}} \ar@{->>}[r]^{R(x)}  & R[f]   \ar@<-1,ex>[d]_{d_0^{f}}\ar@<+1,ex>[d]^{d_1^{f}}\\
   {K[x]\;} \ar@{->>}[d]_{K(f')} \ar[u] \ar@{>->}[r]^{k_{x}} & X'  \ar@{->>}[r]^{x}\ar[u]  \ar@{->>}[d]_{f'} & X  \ar@{->>}[d]^f \ar[u] \\
   {K[y]\;}  \ar@{>->}[r]_{k_{y}}  &    Y'   \ar@{->>}[r]_y   & Y   }
$$
The maps $R(x)$ and $K(f')$ are regular epimorphisms, since the
lower right-hand side square is a regular pushout. The upper row
is left exact by commutation of limits. The map $R(x)$ is
$\Sigma$-special, because $y$ and $x$ are and $\Sigma$ is
point-congruous.

Since $f'$ and $K(f')$ are $\Sigma$-special, the pairs
$(d_0^{f'},s_0^{f'})$ and $(d_0^{K(f')},s_0^{K(f')})$ are in
$\Sigma$. According to the previous proposition $(d_0^{f},
s_0^{f})$ is in $\Sigma$, and hence $f$ is $\Sigma$-special.
\endproof

Now we can state the version for $\Sigma$-protomodular categories
of the so-called \emph{lower $3 \times 3$ lemma}:

\begin{theo} \label{cuicui}
Let $\EE$ be a pointed regular $\Sigma$-protomodular category for
a point-congruous, $2$-regular, equi-consistent class $\Sigma$ of
points. Consider a weakly $\Sigma$-special $3\times 3$ diagram
\eqref{weakly normalized 3x3 diagram}. When the upper row is
$\Sigma$-exact, so is the lower one.
\end{theo}

\proof When $\Sigma$ is point-congruous, the map $K(f) \colon K[x]
\to K[y]$ is $\Sigma$-special, since $f$ and $f'$ are. When the
upper row is left exact, we get $k_{\phi}=k_{K(f)}$. By
Proposition \ref{p4}, the factorization $t$ such that $K(f)=t \,
\phi$ is a monomorphism, and we get $u=k_y \, t$.

Since $f'$ is $\Sigma$-special, and $y$ and $K(x)$ are regular
epimorphisms, Corollary \ref{cor1} implies that the lower
right-hand side square is a regular pushout, so that $K(f)$ is a
regular epimorphism. Accordingly $t$ is a regular epimorphism and
consequently an isomorphism, so that $u$ is the kernel of $y$.

It remains to show that the regular epimorphism $y$ is
$\Sigma$-special. This is a consequence of Corollary \ref{cor2},
since the lower right-hand side square is a regular pushout, $f$
and $f'$ are $\Sigma$-special, and $x$ and $K(x)$ are
$\Sigma$-special, this last point thanks to the assumption that
the upper row is $\Sigma$-exact.
\endproof

There is a last situation dealing with the $3 \times 3$ lemma. We
observed that $\Sigma l/Y$ is finitely complete when $\Sigma$ is
point-congruous, and regular when $\Sigma$ is $2$-regular. From
\cite{Bourn sigma Maltsev}, we know moreover that $\Sigma l/Y$ is
protomodular when $\EE$ is $\Sigma$-protomodular. Recall from
\cite{Bourn 3x3} the following:

\begin{defi}
A category $\EE$ is said to be \emph{quasi-pointed} when it has an
initial object and the terminal map $\tau_0 \colon 0 \into 1$ is a
monomorphism. A category $\EE$ is said to be \emph{sequentiable}
when it is quasi-pointed, protomodular and regular.
\end{defi}

Any fiber $Cat_B$ is quasi-pointed. The category $Grd_B$ of
groupoids with set of objects $B$ is sequentiable. If a category
$\EE$ is pointed regular and $\Sigma$-protomodulaar w.r.t. a
point-congruous and $2$-regular class $\Sigma$, then the category
$\Sigma l/Y$ is sequentiable.

In a sequentiable category, the kernel of a map $f$ is defined as
the pullback of $f \colon X \to Y$ along the monomorphic initial
map $\alpha_Y \colon 0 \into Y$; a sequence is said to be exact
when $f$ is a regular epimorphism and the following pullback is a
pushout as well:
$$\xymatrix{
         {K[f]\;}  \ar@{>->}[r]^{k_f} \ar@{->>}[d] & X  \ar@{->>}[d]^f\\
         {0\;}   \ar@{>->}[r]_{\alpha_Y}  & Y. }
$$
It was proved in \cite{Bourn 3x3} that, in a sequentiable
category, any weakly $3\times 3$ diagram satisfies the $3\times 3$
lemma w.r.t. this extended notion of exact sequence. Whence the
following:

\begin{prop} \label{cuicuicui}
Let $\Sigma$ be a point-congruous and $2$-regular class of points
in a pointed regular $\Sigma$-protomodular category $\EE$. Any
weakly $3\times 3$ diagram \eqref{weakly normalized 3x3 diagram}
in $\EE$ such that the morphisms $y$, $f'$ and $f' \, x = y \, f$
are $\Sigma$-special satisfies the $3\times 3$ lemma.
\end{prop}

\proof The proof is exactly on the same model as the one of
Proposition \ref{p18}. Let us think of diagram \eqref{weakly
normalized 3x3 diagram} as a diagram in the regular category
$\EE/Y'$. In this diagram, from any object there is a unique
morphism to $Y'$, so that any object in $\EE/Y'$ can be identified
with its domain. In our case, any object lies in $\Sigma l/Y'$
except $K[\phi]$ and $U$. When the lower row is exact, $U$ is in
$\Sigma l/Y'$, and so is $K[\phi]$. When the upper row is exact,
the upper left-hand side square is a pullback, since $k_{f'}$ is a
monomorphism. Accordingly, $k_{\phi}$ is also the kernel of $f \,
k_x( = u \, \phi)$, which is $\Sigma$-special, being in $\Sigma
l/Y'$. Hence the factorization $u$ is a monomorphism by
Proposition \ref{p4}, and $R[\phi]=R[\phi \, u]=R[f \, k_x]$ is in
$\Sigma l/Y'$; consequently $U$ is also in $\Sigma l/Y'$ by
Proposition \ref{p3}. In this way, under any of the conditions of
the $3\times 3$ lemma, the whole diagram lies in the sequentiable
category $\Sigma l/Y'$, and the $3\times 3$ lemma holds.
\endproof

\subsection{The fibers $Cat_Y\EE$}

We recalled in Section \ref{section Sigma-Mal'tsev categories}
what is the class $\Sigma$ of points with fibrant splittings in
the category $Cat_Y$. As in $Cat_1=Mon$, there is a description of
these points in terms of Schreier retractions: they are those
points $(F,S) \colon \mathbb A \splito \mathbb B$ which are
equipped with a function $q \colon \mathbb A \to K[F]$ such that
$\phi=q(\phi) \, SF(\phi)$ for all $\phi \in \mathbb B$, and $q(k
\, S(\psi))=k$ for all $(k,\psi)\in K[F] \times \mathbb B $. $q$
is called the \emph{internal Schreier retraction} of the point
$(F,S)$ (notice that consequently any split epimorphism above a
groupoid has cofibrant splittings: in this case we have
$q(\phi)=\phi \, SF(\phi^{-1})$).

\smallskip

Given any category $\EE$, we denote by $(\;)_0 \colon Cat\EE \to
\EE$ the fibration associating with any internal category its
object of objects, and by $Cat_Y\EE$ the fiber above $Y$. When
$\EE$ is a regular category, the category $Cat \EE$ of internal
categories in $\EE$ is no longer regular in general; however so is
any fibre $Cat_Y \EE$, which is also quasi-pointed.

We say that an internal split epimorphic functor $(\underline
F_1,\underline S_1): \underline  X_1\splito \underline Y_1$
\emph{has cofibrant splittings} when there is a morphism $q_1
\colon X_1 \to K[F_1]$ in the underlying category $\EE$ such that
$m_{X_1} \, (k_{F_1} \, q_1,S_1 F_1)=1_{X_1}$ and $q_1 \, m_{X_1}
\, (k_{F_1},S_1)=p_0^{K_{F_1}}$  where $m_{X_1}$ denotes the
composition map of the internal category $\underline X_1$. Denote
by $\Sigma_Y$ the class of points with cofibrant splittings. By
the Yoneda embedding, we get (see \cite{Bourn sigma Maltsev}): 1)
any split epimorphism above an internal groupoid has cofibrant
splittings; 2) the class $\Sigma_Y$ is stable under pullbacks and
point-congruous; and 3) any fiber $Cat_Y\EE$ is
$\Sigma_Y$-protomodular, and accordingly $\Sigma_Y$-Mal'tsev.

\begin{prop}
Given any regular category $\EE$, any fiber $Cat_Y\EE$ is such
that the class of points having cofibrant splittings is
$2$-regular.
\end{prop}

\proof Thanks to the internal Schreier retractions, the same proof
as in Proposition \ref{examp} holds.
\endproof

Accordingly, Theorems \ref{th1} and \ref{th2}, as well as
Proposition \ref{p18}, hold for the fibers $Cat_Y\EE$. As for the
quasi-pointed version of the $3\times 3$ lemma, Propositions
\ref{p22} and \ref{cuicuicui} hold as well. Let us conclude this
section with the following observation:

\begin{prop}
The class $\Sigma_Y$ of points with fibrant splittings, in any
fiber $Cat_Y$, is equi-consistent. Given any category $\EE$, the
same holds for any fiber $Cat_Y\EE$.
\end{prop}

\proof The proof of the first assertion mimics the proof of
Proposition \ref{e-q}. Starting with a parallel pair $(\phi,\phi')
\colon a \rightrightarrows b$ such that $\phi R\phi'$, we have to
show that $q(\phi)Rq(\phi')$. Since $R$ is a
$\Sigma_Y$-equivalence relation, there is an endo map
$\chi(\phi,\phi') \colon b \to b$ such that
$\phi'=\chi(\phi,\phi') \, \phi$ and $1_BR\chi(\phi,\phi')$;
whence $1Rq\chi(\phi,\phi')$, since the point $(K_0 F, K_0 S)$ has
fibrant splittings. As in the proof of Proposition \ref{e-q}, the
identities:
\begin{enumerate}
\item $q(\phi')=q\chi(\phi,\phi') \, q(SF\chi(\phi,\phi') \,
q(\phi))$;

\item $q(SF\chi(\phi,\phi') \, q(\phi)) \,
SF(\chi(\phi,\phi'))=SF\chi(\phi,\phi')$
\end{enumerate}
produce the diagram:
$$\xymatrix{
          q(\phi)   \ar[r]^R \ar@{.>}[rd] & SF\chi(\phi,\phi') \, q(\phi)\\
         q(SF\chi(\phi,\phi') \, q(\phi))  \ar[r]_{R} \ar[ru]^<<<<<<{R}  & q(\phi'), }
       $$
and so $q(\phi)Rq(\phi')$ as desired.

The proof of the second assertion is a consequence of the Yoneda
embedding and of the equations satisfied by the internal Schreier
retractions.
\endproof

Accordingly, Proposition \ref{cuicui} holds for any fibre
$Cat_Y\EE$, when $\EE$ is regular.

\section{A remark on Baer sums} \label{section Baer sums}

In \cite{MMS nine lemma}, given any Schreier exact sequence with
abelian kernel in $Mon$:
$$\xymatrix{
    {A\;}  \ar@{>->}[r]^{k_f}  & X  \ar@{->>}[r]^f & Y }
$$
a monoid action $\phi:Y\to End(A)$ of $Y$ on $A$ is produced,
giving rise to an abelian group $Y\ltimes A \splito Y$ in the
slice category $Mon/Y$. Then a Baer sum construction of exact
sequences  giving rise to the same monoid action as above is
described, together with a so-called push forward construction of
Baer sums.

In \cite{Bourn quandles}, a similar Baer sum construction is given
for quandles, concerning $\Sigma$-special exact sequences with
\emph{abelian kernel equivalence relations}, where $\Sigma$ is the
class of acupuncturing points of quandles. Recall that a
$\Sigma$-special map $f \colon X \to Y$ has an abelian kernel
equivalence relation when there is a Mal'tsev operation $p \colon
R[f] \times_X R[f] \to X$ in $\Sigma l/Y$ (see the proof of the
lemma below for more details).

But, as proved in \cite{MartinsMontoli}, a Schreier exact sequence
in $Mon$ has an abelian kernel if and only if it has an abelian
kernel equivalence relation. We can be even more precise about
this point:

\begin{lma}\label{Schequ}
Given any Schreier equivalence relation $R$ on $X$ in $Mon$, the
following conditions are equivalent:
\begin{itemize}
\item[(i)] $K[d_0^R]$ is an abelian group;

\item[(ii)] for all $1Rt$ and $xRx'$, we get $q(xR(x \cdot
t))=q(x'R(x' \cdot t))$, where $q$ is the Schreier retraction of
the point $(d_0^R, s_0^R)$;

\item[(iii)] the Schreier equivalence relation $R$ is abelian.
\end{itemize}
\end{lma}

\proof We start by observing that the map $q$ satisfies the
following equalities (see Corollary $2.8$ in \cite{MMS Baer
sums}):
\begin{itemize}
\item[(a)] $q(aRb) \cdot b = q(bRa) \cdot b = a$;

\item[(b)] $q((b \cdot b')R(a \cdot a'))=q(bRa) \cdot q(bR(b \cdot
q(b'Ra')))$.
\end{itemize}

Suppose (i). Thanks to the equality (a) above we get:
$$q(xR(x \cdot t)) \cdot x'=q(xR(x \cdot t)) \cdot q(xRx') \cdot x=q(xRx') \cdot q(xR(x \cdot t)) \cdot x=q(xRx') \cdot x \cdot t=x' \cdot t.$$
Whence $q(xR(x \cdot t))=q(x'R(x' \cdot t))$ and (ii).\\

Now suppose (ii). Define the Mal'tsev operation by
$p(aRbRc)=q(bRa) \cdot c$. First we have to show the Mal'tsev
axioms, namely: $p(aRaRc)=c=p(cRaRa)$. We have:
$$q(aRa) \cdot c=1 \cdot c=c=q(aRc) \cdot a=p(cRaRa).$$
It remains to show that $p$ is a monoid homomorphism, namely:
$$q(bRa) \cdot c \cdot q(b'Ra') \cdot c'=q((b \cdot b')R(a \cdot a')) \cdot c \cdot c'.$$ Using the equality (b) above, it is enough
to check the following equality:
$$q(bR(b \cdot q(b'Ra'))) \cdot c=c \cdot q(b'Ra').$$ We have $bRc$, so that, by
(ii), it remains to check that:
$$q(cR(c \cdot q(b'Ra'))) \cdot c=c \cdot q(b'Ra').$$ This is
true since the  point $(d_0^R,s_0^R)$ is a Schreier point. Whence (iii).\\
The fact that condition (iii) implies (i) is true in any category.
\endproof

As shown in \cite{MartinsMontoli}, the fact that a Schreier exact
sequence has an abelian kernel if and only if it has an abelian
kernel equivalence relation is actually true in every category of
monoids with operations in the sense of \cite{MartinsMontoliSobral
mon w op}, in particular in the category $SRng$ of semirings. An
analogue of Lemma \ref{Schequ} can be proved in a
very similar way for all such algebraic structures. \\

According to Lemma \ref{Schequ}, a Schreier exact sequence has an
abelian kernel if and only if it has an  abelian kernel
equivalence relation. So that both Baer sum constructions in $Mon$
and $Qnd$ appear to be of the same nature.

Actually they are both particular instances of a very general
situation described in \cite{Bourn Baer sums} concerning
Barr-exact Mal'tsev categories, which we are now going to recall.
First observe that the category $\Sigma l/Y$, in both cases of
$Mon$ and $Qnd$ is a Barr-exact category by Proposition \ref{p3},
since so are $Mon$ and $Qnd$ because they are both varieties of
universal algebras. Secondly $\Sigma l/Y$, in both cases, is a
Mal'tsev category, since both $Mon$ and $Qnd$ are
$\Sigma$-Mal'tsev ones.

Now, given any Mal'tsev category, we say that an object $X$ is
\emph{affine} when it is endowed with a (necessarily unique)
Mal'tsev operation, namely a ternary operation $p \colon X\times
X\times X \to X$ satisfying $p(x,y,y)=x=p(y,y,x)$. Recall from
\cite{Bourn Baer sums} that a Mal'tsev operation in any Mal'tsev
category is necessarily associative and commutative. Denote by
$Aff\EE$ the full subcategory of affine objects. An object $X$ is
\emph{abelian} when it is endowed with a (necessarily unique and
abelian) internal group structure. An abelian object is nothing
but an affine object $X$ equipped with a point $0_X \colon 1\into
X$; in set-theoretical terms, the abelian group operation $a+b$ is
just $p(a,0,b)$. Denote by $Ab\EE$ the full subcategory of abelian
objects, and by $U \colon Ab\EE \to Aff\EE$ the forgetful functor.
An object $X$ has \emph{global support} if the terminal map
$\tau_X \colon X \to 1$ is a regular epimorphism.

\begin{defi}[\cite{Bourn Baer sums}]
Given a Barr-exact Mal'tsev category $\EE$ and an affine object
$X$ with global support, the \emph{direction} $d(X)$ of $X$ is the
abelian object defined by the following diagram:
  $$ \xymatrix@C=3pc@R=2pc{
    X\times X \times X \ar@<-1,ex>[d]_{p_0}\ar@<+1,ex>[d]^{(p,p_1 p_2)} \ar@<-1,ex>[r]_{(p_0 p_0,p)}\ar@<+1,ex>[r]^{p_2}
    & X\times X \ar@<-1,ex>[d]_{p_0} \ar@<+1,ex>[d]^{p_1} \ar[l] \ar@{.>>}[r]^{q_X} & d(X) \ar@<+1,ex>[d]^{\tau_{d(X)}}\\
    X\times X \ar@<-1,ex>[r]_{p_0} \ar@<+1,ex>[r]^{p_1} \ar[u]_{} & X\ar[u]_{} \ar[l] \ar@{.>>}[r]_{\tau_X} & 1 \ar[u]^{o_X}
}
$$
where $p$ is the Mal'tsev operation and $q_X$ is the quotient of
the upper horizontal equivalence relation.
\end{defi}

The two left-hand side vertical groupoid structures associated
with the structure of equivalence relation give rise (by passage
to the quotient) to a vertical right-hand side groupoid structure
above the terminal object $1$, namely to a group structure in
$\EE$. To make it short, we shall denote by the only symbol $d(X)$
this whole group structure. It is clear that, when $X$ is an
abelian object, we have $d(U(X))\simeq X$ in a natural way.

Take $\EE=\Sigma l/Y$ in $Mon$; the objects of $\EE$ are the
special Schreier homomorphisms $f \colon X \to Y$ with codomain
the monoid $Y$. The objects with global supports are the special
Schreier surjective homomorphisms over $Y$. An object $f$ is
affine in $\EE$ if and only if it has an abelian kernel
equivalence relation, i.e. if and only if is has an abelian
kernel. In this case, its direction is precisely the abelian group
$Y\ltimes A \splito Y$ described at the beginning of this section.
\\

Let us denote by $\EE_g$ the full subcategory of $\EE$ whose
objects are the ones with global support. We have the following:

\begin{prop}[\cite{Bourn Baer sums}, Proposition $6$ and Theorem $7$]
Let $\EE$ be a Barr-exact Mal'tsev category. The construction of
the direction gives rise to a functor \linebreak $d \colon
Aff\EE_g \to Ab\EE$ which preserves finite products and is a
cofibration whose morphisms in the fibers are isomorphims.
\end{prop}

Suppose again that $\EE$ is $\Sigma l/Y$ in $Mon$. The fact that
the morphisms in the fibers are isomorphisms is nothing but the
short five lemma for Schreier exact sequences (see Proposition
$7.2.2$ in \cite{Schreier book} for a proof of this version of the
short five lemma). Moreover, the construction of the cocartesian
map above a morphism in $Ab\EE$ is precisely what is called the
push forward construction described in \cite{MMS nine lemma}.

\begin{prop}[\cite{Bourn Baer sums}, Theorem $9$]
Let $\EE$ be a Barr-exact Mal'tsev category. Given any abelian
object $A$ in $\EE$, the groupoid $d^{-1}(A)$ is canonically
endowed with a closed symmetric monoidal structure, and any
change-of-base functor of the cofibration $d$ is monoidal.
\end{prop}

\proof Given any pair of affine objects $C$ and $C'$, both with
the abelian group $A$ as direction, we first observe that
$d(C\times C')=A\times A$; then their Baer sum is the codomain of
the cocartesian map above the abelian group operation $+ \colon
A\times A \to A$, namely the codomain of the quotient map of the
equivalence relation on $C\times C'$ determined by the following
pullback in $\EE$:
$$ \xymatrix@C=3pc@R=2pc{
    R \ar@{->>}[r]  \ar@{ >->}[d] & A \ar@{ >->}[d]^{(1_A,-1_{A})}\\
    (C\times C')\times (C\times C') \ar@{->>}[r]_-{q_{C\times C'}} & A\times A
}
$$
\endproof

When $\EE$ is the category $\Sigma l/Y$ in $Mon$, the construction
of this tensor product on the fiber $d^{-1}(Y\ltimes A \splito Y)$
coincides with the Baer sum, described in \cite{MMS Baer sums, MMS
nine lemma}, of exact sequences whose associated monoid action
corresponds to the split sequence $Y\ltimes A \splito Y$. It is
the case, as well, of the construction of the Baer sum of special
exact sequences of quandles described in \cite{Bourn quandles}.

\section*{Acknowledgements}

This work was partially supported by the Programma per Giovani
Ricercatori ``Rita Levi Montalcini'', funded by the Italian
government through MIUR.

\end{document}